\documentclass[a4paper,11pt]{amsart}
\usepackage{english}
\usepackage{amsmath,enumerate, amsfonts, amssymb,amsthm, mathptm}
\input xy
\xyoption{all}

\newtheorem{Def}{Definition}[section]
\newtheorem{Theo}[Def]{Theorem}
\newtheorem{Prop}[Def]{Proposition}
\newtheorem{Rem}[Def]{Remark}
\newtheorem{Lem}[Def]{Lemma}
\newtheorem{Cor}[Def]{Corollary}
\newtheorem{Ex}[Def]{Example}

\DeclareMathOperator{\hess}{hess}
\DeclareMathOperator{\pfa}{pf}
\DeclareMathOperator{\tr}{Tr}
\DeclareMathOperator{\dete}{det}

\title{Nondegenerate Monge-Amp\`ere structu\-res in 
  dimension 6}

\author{Bertrand BANOS}
\address{Bertrand Banos\\D\'epartement de Math\'ematiques\\
  Universit\'e d'Angers\\ 2 bd La\-voi\-sier, 49045 Angers, France}
\email{bertrand.banos@univ-angers.fr}

\begin{document}

\begin{abstract}
We define a nondegenerate Monge-Amp\`ere structure on a
$6$-di\-men\-sio\-nal mani\-fold as a pair $(\Omega,\omega)$, such that
$\Omega$ is a symplectic form and $\omega$ is a  $3$-differential form
which satisfies $\omega\wedge\Omega=0$ and which is nondegenerate in
the sense of Hitchin. We associate with such a pair a generalized
almost (pseudo) Calabi-Yau structure and we study its integrability
from the point of view of Monge-Amp\`ere operators theory. The result
we prove appears as an  analogue of Lychagin and Roubtsov theorem on
integrability of the almost complex or almost product structure
associated with an elliptic or hyperbolic Monge-Amp\`ere equation in
the dimension $4$. We study from this point of view the example of the
Stenzel metric on $T^*S^3$.  

\emph{Key words}: Calabi-Yau manifolds, special lagrangian
submanifolds, Monge-Amp\`ere equations, Symplectic forms

\emph{AMS classification}: 34A26, 58A10, 53D05, 32Q60, 14J32   

\end{abstract}

\maketitle

\section{Introduction}

A Monge-Amp\`ere equation is a differential equation which is
nonlinear in a very specific way: its nonlinearity is the determinant
one. In the dimension $2$, such an equation can be written as
\begin{equation}{\label{EMA}}
A\frac{\partial^2 f}{\partial q_1^2}+2B\frac{\partial^2 f}{\partial 
  q_1\partial q_2}+ C\frac{\partial^2 f}{\partial q_2^2}+
D\Big(\frac{\partial^2 f}{\partial q_1^2}\frac{\partial^2 f}{\partial
  q_2^2}-(\frac{\partial^2 f}{\partial q_1\partial q_2})^2\Big)+E=0, 
\end{equation}
where $A,B,C,D$ and $E$ are smooth functions on the jet space
$J^1\mathbb{R}^2$. This equation is said to be symplectic if these
coefficients are actually functions on the quotient bundle
$J^1\mathbb{R}^2/J^0\mathbb{R}^2$, which can be identified  with the
cotangent bundle $T^*\mathbb{R}^2$.      

The Monge-Amp\`ere operators theory proposed by Lychagin (\cite{L})
associates to each symplectic Monge-Amp\`ere equation on an
$n$-dimensional manifold $M$ a pair of differential forms
$(\Omega,\omega)$ on the cotangent bundle $T^*M$ , where $\Omega$ is
the canonical symplectic form and $\omega\in \Omega^n(T^*M)$ is an
$\Omega$-effective form, i.e. $\omega\wedge\Omega=0$. To be more precise,
the symplectic Monge-Amp\`ere equation associated with such a pair
$(\Omega,\omega)$ is the differential equation
\begin{equation}{\label{MAO}}
(df)^*(\omega)=0,
\end{equation}
where $df:M\rightarrow T^*M$ is the natural section defined by a
smooth function $f$ on $M$. For instance, \eqref{EMA} is 
associated with the symplectic form on $T^*\mathbb{R}^2$
$$
\Omega= dq_1\wedge dp_1+dq_2\wedge dp_2,
$$
and the differential form
$$
\begin{aligned}
\omega&=A dp_1\wedge dq_2 + B(dq_1\wedge dp_1-dq_2\wedge dp_2)+C
dq_1\wedge dp_2\\ 
&+D dp_1\wedge dp_2+ E dq_1\wedge dq_2.\\ 
\end{aligned}
$$   

A generalized solution of the differential equation \eqref{MAO} is a
lagrangian submanifold $L$ of $(T^*M,\Omega)$ on which $\omega$
vanishes. $L$ is locally the graph of a section
$df:M\rightarrow T^*M$ with $f$ solution of $\eqref{MAO}$ if and only
if the projection $L\rightarrow M$ is locally a diffeomorphism.   

Generalizing this notion, we define a symplectic Monge-Amp\`ere
structure on a $2n$-dimen\-sional $X$ as a pair of differential forms
$(\Omega,\omega)$, $\Omega\in \Omega^2(X)$ being symplectic and
$\omega\in \Omega^n(X)$ being $\Omega$-effective. In the dimension
$n=2$, if $\omega$ is nondegenerate (i.e. the pfaffian
$\pfa(\omega)=\frac{\omega\wedge\omega}{\Omega\wedge\Omega}$ is non
zero), the equality
$$
\frac{\omega}{\sqrt{|\pfa(\omega)|}}=\Omega(A_\omega .,.)
$$
defines a section $A_\omega: X\rightarrow TX\otimes T^*X$ which is
either an almost complex structure ($A_\omega^2=-Id$) or an almost
product structure ($A_\omega^2=Id$) on $X$. Lychagin and Roubtsov gave
in \cite{LR1} a necessary and sufficient condition for $A_\omega$
to be integrable:

\begin{Prop}{\label{int2}}
$A_\omega$ is integrable if and only if 
$$
d\Big(\frac{\omega}{\sqrt{|\pfa(\omega)|}}\Big)=0.
$$
\end{Prop}

From the point of view of differential equation (i.e. the local point
of view), this result can be formulated as follows:

\begin{Prop}{\label{class2}}
A symplectic Monge-Amp\`ere equation $\Delta_\omega=0$ on
$\mathbb{R}^2$ is symplectically equivalent to one of these two
equations  
$$
\begin{cases}
\Delta f=0,&(\pfa(\omega)>0)\\
\square f=0,&(\pfa(\omega)<0)\\
\end{cases}
$$
if and only if 
$$
d\Big(\frac{\omega}{\sqrt{|\pfa(\omega)|}}\Big)=0.
$$
\end{Prop}

To classify the symplectic Monge-Amp\`ere equations with constant
coefficients in the dimension $3$, Lychagin and Roubsov introduced a
quadratic invariant $q_\omega$ associated with each effective $3$-form
$\omega$ (\cite{LR3}). Hitchin has defined a linear invariant
$K_\omega$ associated with each $3$-form $\omega$ (\cite{Hi}) in order
to study the geometry of $3$-dimensional Calabi-Yau manifolds. We show
that these two invariants coincide in the effective case
(proposition \ref{compatibilite}) . This remark allows us to combine
Lychagin, Roubtsov and Hitchin works to demonstrate a result analogous
to \ref{int2} and \ref{class2} in the dimension $3$. Let us introduce
some notations to state this result. We suppose that our $3$-form
$\omega\in \Omega^3(X^6)$ is nondegenerate in the Hitchin sense. The
Hitchin pfaffian $\lambda(\omega)$ is then nonzero and $\omega$ can be
decomposed in an unique way as a sum of two complex decomposable
$3$-forms: $\omega=\alpha+\beta$. Following Hitchin, we denote by
$\hat{\omega}$ the dual form associated with $\omega$. We associate to
the pair $(\Omega,\omega)\in \Omega^2(X^6)\times\Omega^3(X^6)$ a
geometric structure on $X$ which we call \emph{generalized almost
  Calabi-Yau structure}. This structure is essentially composed of a
(pseudo) metric $q_\omega$, an almost complex structure or almost
product structure $K_\omega$ which is compatible with $q_\omega$ and
$\Omega$ and two decomposable $3$-forms whose associated distributions
are distributions of $K_\omega$ eigenvectors. We extend Hitchin's
results to demonstrate the analogous of \ref{int2}: 

\begin{Prop}{\label{int3}}
The generalized almost Calabi-Yau structure 
$$
(q_\omega,K_\omega,\Omega,\alpha,\beta)
$$
is ``integrable'' if and
only if 
$$
\begin{cases}
d\Big( \frac{\omega}{\sqrt[4]{|\lambda(\omega)|}}\Big)=0,&\\
d\Big( \frac{\hat{\omega}}{\sqrt[4]{|\lambda(\omega)|}}\Big)=0.&\\ 
\end{cases}
$$
\end{Prop}

We demonstrate then the local version of this proposition:

\begin{Theo}{\label{theo}}
A symplectic Monge-Amp\`ere equation in the dimension $3$ associated
with a nondegenerate Monge-Amp\`ere structure $(\Omega,\omega)$ is
symplectically equivalent to one of these three equations
$$
\begin{cases}
\hess(f)=1,&\\
\Delta f-\hess(f)=0,&\\
\square f +\hess(f)=0,&\\
\end{cases}
$$
if and only if
$$
\begin{cases}
d\Big( \frac{\omega}{\sqrt[4]{|\lambda(\omega)|}}\Big)=0,&\\
d\Big( \frac{\hat{\omega}}{\sqrt[4]{|\lambda(\omega)|}}\Big)=0,&\\ 
q_\omega\text{ is flat}.&\\
\end{cases}
$$
\end{Theo}

Our motivation is to generalize the notion of ``special lagrangian
submanifolds'' in the dimension $3$. These were first introduced
by Harvey and Lawson in their famous paper \emph{Calibrated
  Geometries} (\cite{HL}) as examples of minimal submanifolds. Recall
that a $p$-calibration on a riemannian manifold $(Y,g)$ is a closed
differential $p$-form $\phi\in \Omega^p(Y)$ such that for any point
$y\in Y$ and any oriented $p$-plane $V$ of $T_yY$, the following
inequality
$$
\phi_y|_V\leq vol_V
$$
holds. Here $vol_V$ is the volume exterior form on $V$ defined by the
metric $g$ and the orientation on $V$. An oriented  $p$-dimensional
submanifold $L$ is said to be $\phi$-calibrated if for any $y\in L$ 
$$
\phi_y|_{T_yL}=vol_{T_yL}.
$$
Calibrated submanifolds are volume-minimizing in their homology
classes. The real form $Re(\alpha)$ with $\alpha=dz_1\wedge\ldots\wedge
dz_n$ is an example of $n$-calibration on $\mathbb{C}^n$ and
$Re(\alpha)$-calibrated sub\-ma\-ni\-folds are said to be special
lagrangian. This notion of special lagrangian calibration can be
generalized on Calabi-Yau manifolds, i.e. K\"ahler manifolds
endowed with an holomorphic volume form $\alpha\in
\Omega^{n,0}(Y)$ such that
$$
\frac{\alpha\wedge\overline{\alpha}}{\Omega^n}\text{ is constant},
$$     
$\Omega$ being the K\"ahler form. These special lagrangian
submanifolds attracted a lot of  attention of many mathematicians
in the last few years after Strominger, Yau and Zaslow proposed a
geometric construction of mirror manifolds based on the conjecture of
existence  of toric special lagrangian fibration (\cite{SYZ}).  

Gromov noted in a discussion with Roubtsov that the Monge-Amp\`ere
structures can be seen as an analogue of the calibrations. Effective forms
correspond to calibrations and lagrangian submanifolds correspond to
calibrated submanifolds. Moreover, special lagrangian submanifolds are
in the intersection of these two approach. Harvey and Lawson  have
shown that special lagrangian submanifolds of
$(\mathbb{C}^n,Re(\alpha))$ are actually the generalized solutions of the
differential equation associated with the Monge-Amp\`ere structure
$(\Omega,Im(\alpha))$ with 
$$
\Omega=\frac{i}{2}(dz_1\wedge
d\overline{z_1}+\ldots+dz_n\wedge d\overline{z_n}). 
$$
In the dimension $3$ this equation is
\begin{equation}{\label{slag}} 
\Delta f-\hess (f)=0.
\end{equation}

Our aim is to show that the problem of local equivalence of
Monge-Amp\`ere equations in the dimension $3$ is the local expression
of the problem of integrability of some geometrical structures that
generalize in a very natural way the Calabi-Yau structure. This
approach gives a different  description of  Calabi-Yau manifolds,
seeing them more as symplectic manifolds than as complex manifolds. 

In the first section we recall the Lychagin's approach to
Monge-Amp\`ere equations. We study as an example Chynoweth-Sewell's
equations which come from the ``semi-geostrophic'' model of
athmosphere dynamics. We remark that
they are all equivalent to the classic Monge-Amp\`ere equation
$\hess(f)=1$.
In the second section we adapt Hitchin's works on $3$-forms to the
effective case. 
In the last section, we define the notion of generalized almost
Calabi-Yau structure and  demonstrate the theorem \ref{theo}. We study
as an example the Stenzel metric on $T^*S^3$.

This article constitutes a part of the author's PhD thesis being
prepared at Angers University. I would like to thank my advisor
Volodya Roubtsov for suggesting the problem and helpfull discussions. I
would like also to thank Oleg Lisovyy for all his help. I am very
grateful to professor Benjamin Enriquez for his valuable remarks and
suggestions. I would like finally to thank Mich\`ele Audin for the
lecture she gave in Barcelona summer school (\cite{A}) which helped me to
understand the geometry of lagrangian and special lagrangian
submanifolds.     

\section{Effective forms and Monge-Amp\`ere operators}

Let $(V,\Omega)$ be a symplectic $2n$-dimensional vector space over
$\mathbb{R}$ and $\Lambda^*(V^*)$ the space of exteriors forms on
$V$. Let  $\Gamma:V\rightarrow V^*$ be the isomorphism
determined by $\Omega$ and let $X_\Omega\in\Lambda^2(V)$ be the unique
bivector such that $\Gamma^*(X_\Omega)= \Omega$.

Following Lychagin (see \cite{L}), we introduce the operators $\top:
\Lambda^k(V^*)\rightarrow \Lambda^{k+2}(V^*)$, $\omega\mapsto
\omega\wedge \Omega$ and  $\bot:\Lambda^k(V^*)\rightarrow
\Lambda^{k-2}(V^*)$, $\omega\mapsto i_{X_\Omega}(\omega)$. They
have the followings properties:
$$
\begin{cases}
[\bot,\top](\omega)=(n-k)\omega \;\text{, $\forall \omega\in
  \Lambda^k(V^*)$};&\\ 
\bot: \Lambda^k(V^*)\rightarrow \Lambda^{k-2}(V^*)\text{ is into for }
k\geq n+1;&\\ 
\top: \Lambda^k(V^*)\rightarrow \Lambda^{k+2}(V^*)\text{ is into for }
k\leq n-1.&\\ 
\end{cases}
$$
We will say that a $k$-form $\omega$ is effective if $\bot\omega=0$
and we will denote by $\Lambda^{k}_\varepsilon(V^*)$ the vector space
of effective $k$-forms on $V$. When  $k=n$, $\omega$ is effective if
and only if $\omega\wedge \Omega=0$. 

The next theorem explains the fundamental role played by the effective
forms in the theory of Monge-Amp\`ere operators (see \cite{L}):

\begin{Theo}[Hodge-Lepage-Lychagin]{\label{hodge}}
\begin{enumerate}
\item Every form $\omega\in \Lambda^k(V^*)$ can be u\-ni\-que\-ly decomposed
  into the finite sum
$$
\omega= \omega_0+\top \omega_1+\top^2\omega_2+\ldots,
$$
where all $\omega_i$ are effective forms.
\item If two  effective $k$-forms vanish on the same $k$-dimensional
  isotropic vector subspaces in $(V,\Omega)$, they are proportional. 
\end{enumerate}
\end{Theo}

Let $M$ be an $n$-dimensional smooth manifold. Denote by $J^1M$ the
space of $1$-jets of smooth functions on $M$ and by $j^1(f):
M\rightarrow J^1M$,  $x\mapsto [f]_{x}^1$ the natural section
associated with a smooth function $f$ on $M$. The  Monge-Amp\`ere
operator 
$$
\Delta_\omega: C^\infty(M)\rightarrow \Omega^n(M)
$$
associated with a differential $n$-form $\omega\in \Omega^n(J^1M)$ is
the differential operator  
$$
\Delta_\omega(f)=j_1(f)^*(\omega).
$$

Let $U$ be the contact $1$-form on $J^1M$ and $X_1$ be the Reeb's vector
field. Denote by $C(x)$ the kernel of  $U_x$ for $x\in
J^1M$. $(C(x),dU_x)$ is a $2n$-dimensional symplectic vector space and
$$
T_xJ^1M= C(x)\oplus \mathbb{R}X_{1x}.
$$
A generalized solution of the equation $\Delta_\omega=0$ is a
legendrian submanifold $L^{n}$ of $(J^1M,U)$ such that
$\omega|_L=0$. Note that  $T_xL$ is a
lagrangian subspace of $(C(x),dU_x)$ in each point $x\in L$, and that  $L$ is
locally the graph of a section $j^1(f)$, where $f$ is a regular
solution of the equation $\Delta_\omega(f)=0$, if and only if the projection $\pi:J^1M\rightarrow M$ is a local diffeomorphism on $L$. 

We will denote by $\Omega^*(C^*)$ the space of differential forms
vanishing on $X_1$. In each point $x$, $(\Omega^k(C^*))_x$ can be
naturally identified with $\Lambda^k(C(x)^*)$.  Let
$\Omega_\varepsilon^*(C^*)$ be the space of forms which are effective
on $(C(x), dU_x)$ in each point $x\in J^1M$. The first part of the
theorem \ref{hodge} means that  
$$
\Omega^*_\varepsilon(C^*)=\Omega^*(J^1M)/I_C,
$$
where $I_C$ is the Cartan ideal generated by $U$ and $dU$. The second
part means that two differential $n$-forms $\omega$ and $\theta$ on
$J^1M$ determine the same Monge-Amp\`ere operator if and only if
$\omega-\theta\in I_C$.

$Ct(M)$, the pseudo-group of contact diffeomorphisms on $J^1M$,
naturally acts on the set of Monge-Amp\`ere operators in the following
way 
$$
F(\Delta_\omega)=  \Delta_{F^*(\omega)},
$$
and the corresponding infinitesimal action is
$$
X(\Delta_\omega) = \Delta_{L_X(\omega)}.
$$
We are interested in a more restrictive class of operators, the class
of  symplectic operators. These operators satisfy  
$$
X_1(\Delta_\omega) = \Delta_{L_{X_1}(\omega)}=0.
$$
Let $T^*M$ be the cotangent space and $\Omega$ be the canonical
symplectic form on it. Let us consider the projection $\beta:
J^1M\rightarrow T^*M$, defined by the following commutative diagram: 
$$
\xymatrix{
  \mathbb{R}&&J^1M\ar[ll]_\alpha\ar[rr]^\beta&&T^*M\\
&&&&\\
&&M\ar[lluu]^f\ar[uu]_{j^1(f)}\ar[rruu]_{df}&&\\
}
$$

We can naturally identify the space $\{\omega\in
\Omega^*_\varepsilon(C^*): L_{X_1}\omega=0\}$ with the space of
effective forms on $(T^*M,\Omega)$ using this projection $\beta$. Then,
the group acting on these forms is the group of symplectomorphisms of
$T^*M$. 

\begin{Def}
A Monge-Amp\`ere structure on a $2n$-dimensional manifold $X$ is a
pair of differential form $(\Omega,\omega)\in
\Omega^2(X)\times\Omega^n(X)$ such that $\Omega$ is symplectic and
$\omega$ is $\Omega$-effective i.e. $\Omega\wedge\omega=0$.
\end{Def}

When we locally identify the symplectic manifold  $(X,\Omega)$ with
$(T^*\mathbb{R}^n,\Omega_0)$, we can then associate to the pair
$(\Omega,\omega)$ a symplectic Monge-Amp\`ere equation
$\Delta_\omega=0$. Conversely, any symplectic Monge-Ampere
equation $\Delta_\omega=0$ on a manifold $M$ is associated with 
Monge-Amp\`ere structure $(\Omega,\omega)$ on $T^*M$.

\begin{Ex}  The Chynoweth-Sewell's equations are an example of
  Monge-Amp\`ere equations with constant coefficients. They come from
  the ``semi-geo\-stro\-phic  model'' of Atmosphere Dynamics (\cite{CS}):
\begin{equation}{\label{chy-sew}}
\frac{\partial^2 f}{\partial x^2}\frac{\partial^2 f}{\partial
  y^2}-(\frac{\partial^2 f}{\partial x\partial y})^2+\frac{\partial^2
  f}{\partial z^2}=\gamma, \;\; \gamma\in \mathbb{R}
\end{equation} 
Note that $F(x,y)-\frac{1}{2}z^2$ is a solution of $\eqref{chy-sew}$
for the particular case $\gamma=0$ when $F$ is a solution of
$\hess(F)=1$. For instance,
$$
\frac{1}{3}\sqrt{(x^2+2y)^3}-\frac{1}{2}z^2 
$$
is a solution of $\eqref{chy-sew}$ when $\gamma=0$.

The effective form associated with $\eqref{chy-sew}$ is 
$$
\omega= dp\wedge dq\wedge dz+dx\wedge dy\wedge dh-\gamma dx\wedge
dy\wedge dz, 
$$ 
where $(x,y,z,p,q,h)$ is the canonical coordinates system of
$T^*\mathbb{R}^3$. This form is clearly the sum of two decomposable
$3$-forms:
$$
\omega=dp\wedge dq\wedge dz+dx\wedge dy\wedge(dh-\gamma dz).
$$
Then $\phi^*(\omega)=dp\wedge dq\wedge dh -dx\wedge dy \wedge dz$
where $\phi$ is the symplectomorphism
$$
\phi(x,y,z,p,q,h)=(x,y,h,p,q,\gamma h- z).
$$
In other words, Chynoweth-Sewell's equations are symplectically
equivalent to the equation 
\begin{equation}{\label{hess}}
\hess(f)=1.
\end{equation}
It is easy to check that 
$$
f(x,y,z)=\int_{a}^{\sqrt{xy+yz+zx}} (b+4\xi^3)^{1/3}d\xi 
$$
is a regular solution of $\eqref{hess}$. Therefore, 
$$
L=\Big\{(x,y,(x+y)\alpha,(y+z)\alpha,(z+x)\alpha,\gamma(x+y)\alpha-z)\Big\}
$$
is an example of generalized solution of $\eqref{chy-sew}$ with 
$$
\alpha=\frac{1}{2}(\frac{b}{(xy+yz+zx)^{\frac{3}{2}}}+4)^{\frac{1}{3}}.
$$
\end{Ex}

\section{The geometry of effective $3$-forms in the dimension $6$}

\subsection{The action of $SL(6,\mathbb{R})$ on
  $\Lambda^3(\mathbb{R}^6)$} 

We recall first Hitchin's results on the geometry of $3$-forms. Let
$V$ be a $6$-dimensional real vector space.  We denote by
$A:\Lambda^5(V^*)\rightarrow V\otimes \Lambda^6(V^*)$ the isomorphism
induced by the exterior product and we fix a volume form $\theta$
on $V$. The linear map $K_\omega^\theta: V\rightarrow V$ associated with
$\omega\in\Lambda^3(V^*)$ is defined by
$$
K_\omega^\theta(X)\theta = A(i_X(\omega)\wedge\omega).
$$   

\begin{Def}
The Hitchin pfaffian of a $3$-form $\omega\in\Lambda^3(V^*)$ is 
$$
\lambda_\theta(\omega)=\frac{1}{6}\tr(K_\omega^\theta\circ
K_\omega^\theta).
$$
If $\lambda_\theta(\omega)$ is nonzero then $\omega$ is said to be
nondegenerate. 
\end{Def}

\begin{Prop}[Hitchin]
Let $\omega\in \Lambda^3(V^*)$  be nondegenerate. Then,
\begin{enumerate}
\item $K_\omega^\theta\circ
K_\omega^\theta=\lambda_\theta(\omega)Id$.
\item $\lambda_\theta(\omega)>0$ if and only if $\omega=\alpha+\beta$
  where $\alpha$ and $\beta$ are real decomposable $3$-forms on
  $V$. Moreover, if we impose $\frac{\alpha\wedge\beta}{\theta}>0$
  then $\alpha$ and $\beta$ are unique:
$$
\begin{cases}
2\alpha=\omega+
|\lambda_\theta(\omega)|^{-\frac{3}{2}}(K_\omega^\theta)^*(\omega)&\\ 
2\beta=\omega-
|\lambda_\theta(\omega)|^{-\frac{3}{2}}(K_\omega^\theta)^*(\omega)&\\ 
\end{cases}
$$
\item $\lambda_\theta(\omega)<0$ if and only if
  $\omega=\alpha+\overline{\alpha}$ where $\alpha$ is a
  complexe decomposable $3$-form on $V$. Moreover, if we impose
  $\frac{\alpha\wedge\overline{\alpha}}{i\theta}>0$ then $\alpha$ is
  unique: 
$$
\alpha=
\omega+i|\lambda_\theta(\omega)|^{-\frac{3}{2}}K_\omega^*(\omega). 
$$
\end{enumerate}
\end{Prop}

\begin{Rem}
Let $(e_1,\ldots,e_6)$ be a basis of $V$ and fix
$\theta=e_1^*\wedge\ldots \wedge e_6^*$.
\begin{enumerate}
\item $\lambda_\theta(\omega)>0$ if and only if $\omega$ is in the
  $GL(6)$-orbit of
$$
e_1^*\wedge e_2^*\wedge e_3^*+ e_4^*\wedge e_5^*\wedge e_6^*. 
$$
\item $\lambda_\theta(\omega)<0$ if and only if $\omega$ is in the
  $GL(6)$-orbit of 
$$
(e_1^*+ie_4^*)\wedge (e_2^*+ie_5^*)\wedge
(e_3^*+ie_6^*)+(e_1^*-ie_4^*)\wedge (e_2^*-ie_5^*)\wedge (e_3^*-ie_6^*).
$$
\end{enumerate}
Therefore, the action of $GL(6)$ on $\Lambda^3(V^*)$ has two open
orbits separated by the quartic hypersurface $\lambda_\theta=0$. This
explains this notion of nondegenerate $3$-form.
\end{Rem}

\begin{Def}[Hitchin]
Let $\omega$ be a nondegenerate $3$-form on $V$. The dual form
$\hat{\omega}$ is
\begin{enumerate}
\item $\hat{\omega}=\alpha-\beta$ if $\omega=\alpha+\beta$,
\item $\hat{\omega}=i(\overline{\alpha}-\alpha)$ if
  $\omega=\alpha+\overline{\alpha}$. 
\end{enumerate}
\end{Def}

To conclude we remark that the exterior product defines a  symplectic
form on $\Lambda^3(V^*)$ 
$$
\Theta_\theta(\omega,\omega')=\frac{\omega\wedge\omega'}{\theta}
$$
and that the action of $SL(6)$ is hamiltonian:

\begin{Prop}[Hitchin]
The action of $SL(6)$ on $(\Lambda^3(V^*),\Theta_\theta)$ is
hamiltonian with moment map $K^\theta: \Lambda^3(V^*)\rightarrow sl(6)$.
\end{Prop}

\begin{Ex} This invariant $K$ can be used to construct some almost complex
  or almost product structures  on $6$-dimensional
  manifolds. We study in this example the restriction of the famous
  associative $3$-form to the sphere $S^6$. Let $(\mathbb{O},<,>)$ be
  the octonions normed algebra and denote by $E_7$ the $7$-dimensional subspace of imaginary
octonions. The associative form $\phi\in \Lambda^3(E_7^*)$ is defined
by
$$
\phi(x,y,z)=<x,yz>.
$$
(see for instance for instance \cite{HL}). Let us see the sphere $S^6$
as a submanifold of $E_7$ and let us  consider the form $\omega\in
\Omega^3(S^6)$ defined by 
$$
\omega=\frac{1}{\sqrt{2}} \phi|_{S^6}.
$$
Let $\theta$ be the volume form defined by the metric induced on
$S^6$. A straightforward computation shows that
$$
\lambda_\theta(\omega)=-1.
$$
Therefore $K_\omega^\theta$ is an almost complex structure on
$S^6$. This almost complex structure is in fact already known. It
actually coincides with the almost complex structure 
$$
I_x(Y)=x.Y,
$$
with $x\in S^6$ and $Y\in T_xS^6=\big\{Y\in E_7: <x,Y>=0\big\}$. This
almost complex structure is not integrable (\cite{HC}).        
\end{Ex}

\subsection{The action of $SP(3)$ on $\Lambda^3_\varepsilon(V^*)$}

We assume now that $V$ is a $6$-dimensional symplectic vector
space. We fix $\theta=-\frac{1}{6}\Omega^3$ with $\Omega$ the
symplectic form on $V$. We denote $\lambda=\lambda_\theta$,
$K=K^\theta$ and $\Theta=\Theta_\theta$. $\omega\in \Lambda^3(V^*)$ is
said to be effective if $\Omega\wedge\omega=0$. We denote by
$\Lambda^3_\varepsilon(V^*)$ the space of effective
$3$-forms. $(\Lambda_\varepsilon^3(V^*),\Theta)$ is a symplectic
subspace of  $(\Lambda^3(V^*),\Theta)$ since, according to the
Hodge-Lepage-Lychagin theorem, any $3$-form $\omega$ admits the
decomposition  
$$
\omega=\omega_0+\Omega\wedge\omega_1,
$$
with $\omega_0$ effective.

Denote by $sp(3)$ the Lie algebra of $SP(3)=SP(\Omega)$. A
straightforward comuptation shows the following lemma.

\begin{Lem}
Let $\omega$ be a $3$-form on $V$. $\omega$ is effective if and only
if $K(\omega)\in sp(3)$.
\end{Lem}

\begin{Cor}
  The action of $SP(3)$ on $(\Lambda_\varepsilon^3(V^*),\Theta)$ is
  hamiltonian with moment map $K:
  \Lambda_\varepsilon^3(V^*)\rightarrow sp(3)$.   
\end{Cor}
    
Lychagin and Roubtsov have defined an other invariant
$q_\omega\in S^2(V^*)$ associated with one effective $3$-form $\omega$
(\cite{LR3}). It is natural to ask what is the link between these two
invariants. The quadratic form $q_\omega$  is defined by 
$$
q_\omega(X)=-\frac{1}{4}\bot^2(i_X\omega\wedge i_X\omega).
$$
In fact, this invariant gives us the roots of the characteristic
polynom of $i_X\omega$:
$$
(i_X\omega -\xi\Omega)^3=
-\xi(\xi-\sqrt{q_\omega(X)})(\xi+\sqrt{q_\omega(X)})\Omega^3. 
$$
Using this invariant, Lychagin and Roubtsov have listed the different
orbits of the action of $SP(3)$ on $\Lambda^3_\varepsilon(V^*)$. This list has
been completed by the author in \cite{Ba1} and is summed up in table
\ref{table1}.

Computing $q_\omega$ and $K_\omega$ for each normal form and using
their invariant pro\-perties, we can  check the following: 

\begin{Prop}{\label{compatibilite}}
Let $\omega$ be an effective $3$-form on $V$. Then 
$$
q_\omega(X)=\Omega(K_\omega X,X),
$$
for all $X\in V$.
\end{Prop}

\begin{Rem}
The Lie algebra $(sp(3),[\;,\;])$ can be naturally identified with the
Lie algebra $(S^2(V^*),\{\;,\;\})$ where $\{\;,\;\}$ is the Poisson
bracket associated with $\Omega$. $q:
\Lambda^3_\varepsilon(V^*)\rightarrow S^2(V^*)$ is then the moment map
of the hamiltonian action of $SP(3)$ on
$(\Lambda^3_\varepsilon(V^*),\Theta)$. 
\end{Rem}

\begin{table}[!ht]
\begin{small}
\begin{tabular}{|c|c|c|c|c|}
\hline
&$\Delta_\omega=0$ &$\omega$ & $\varepsilon(q_\omega)$ &
$\lambda(\omega)$\\ 
\hline
 1& $\hess(f)=1$ & $dq_{1}\wedge dq_{2}\wedge dq_{3}+ \gamma
dp_{1}\wedge dp_{2}\wedge dp_{3}$&$(3,3)$&$\gamma^4$\\
\hline
 2& $\Delta f- \hess(f)=0$ & $dp_{1}\wedge dq_{2}\wedge dq_{3}- dp_{2}\wedge
dq_{1}\wedge
dq_{3}$&$(0,6)$&$-\gamma^4$\\
&&$+ dp_{3}\wedge dq_{1}\wedge dq_{2}- \nu^2 dp_{1}\wedge
dp_{2}\wedge dp_{3}$&&\\
\hline 
3& $\square f +\hess(f)=0$ & $dp_{1}\wedge dq_{2}\wedge dq_{3}+ dp_{2}\wedge
dq_{1}\wedge
dq_{3}$&$(4,2)$&$-\gamma^4$\\
&&$+ dp_{3}\wedge dq_{1}\wedge dq_{2}+ \nu^2 dp_{1}\wedge
dp_{2}\wedge dp_{3}$&&\\
\hline 
4& $\Delta f=0$ & $dp_{1}\wedge dq_{2}\wedge dq_{3}- dp_{2}\wedge
dq_{1}\wedge dq_{3}+ dp_{3}\wedge dq_{1}\wedge dq_{2}$&
$(0,3)$&$0$\\
\hline
5&$\square f=0$ &$dp_{1}\wedge dq_{2}\wedge dq_{3}+ dp_{2}\wedge
dq_{1}\wedge dq_{3}+ dp_{3}\wedge dq_{1}\wedge dq_{2}$&
$(2,1)$&$0$\\
\hline
6&$\Delta_{q_2,q_3} f=0$&$dp_{3}\wedge dq_{1}\wedge dq_{2}- dp_{2}\wedge
dq_{1}\wedge
dq_{3}$&$(0,1)$&$0$\\
\hline
7&$\square_{q_1,q_2} f=0$&$dp_{3}\wedge dq_{1}\wedge dq_{2}+ dp_{2}\wedge
dq_{1}\wedge
dq_{3}$&$(1,0)$&$0$\\
\hline
8&$\frac{\partial^2 f}{\partial q_1^2}=0$&$dp_{1}\wedge dq_{2}\wedge dq_{3}$&$(0,0)$&$0$\\
\hline
9&& $0$&$(0,0)$&$0$\\
\hline
\end{tabular}
\end{small}
\caption{Classification of effective $3$-forms in the dimension
  $6$}\label{table1}
\end{table}

\section{Geometrical structures associated with nondegenerate
  Monge-Am\-p\`e\-re equations}

\begin{Def}
A Monge-Amp\`ere structure $(\Omega,\omega)$ on a $6$-dimensional
manifold $X$ is called
\begin{enumerate}
\item nondegenerate if $\lambda(\omega)$ never vanishes,
\item elliptic if $\lambda(\omega)<0$ everywhere,
\item hyperbolic if $\lambda(\omega)>0$ everywhere.
\end{enumerate}
A nondegenerate Monge-Amp\`ere structure $(\Omega,\omega)$ is said to be
\begin{enumerate}
\item closed if 
$$
\begin{cases}
d\Big( \frac{\omega}{\sqrt[4]{|\lambda(\omega)|}}\Big)=0,&\\
d\Big( \frac{\hat{\omega}}{\sqrt[4]{|\lambda(\omega)|}}\Big)=0.&\\ 
\end{cases}
$$
\item locally constant if there exists a Darboux coordinates system of
  $(X,\Omega)$ in which $\omega$ has constant coefficients.
\end{enumerate}
\end{Def}

\subsection{Generalized Calabi-Yau structures}

\begin{Def}
A generalized almost Calabi-Yau struc\-tu\-re on a $6$-di\-men\-sional
manifold $X$ is a $5$-uple $(g,\Omega,K,\alpha,\beta)$ where  
\begin{enumerate}
\item $g$ is a (pseudo) metric on $X$,
\item $\Omega$ is a symplectic on $X$,
\item $K$ is a smooth section $X\rightarrow TX\otimes T^*X$ such that
  $K^2=\pm Id$ and such that 
$$
g(U,V)=\Omega(KU,V)
$$
for all tangent vectors $U,V$,
\item $\alpha$ and $\beta$ are (eventually complex) decomposable $3$-forms
   whose associated distributions are the
  distributions of $K$ eigenvectors and such that
$$
\frac{\alpha\wedge\beta}{\Omega^3}\text{ is constant}.
$$
\end{enumerate}
A generalized Calabi-Yau structure $(g,\Omega,K,\alpha,\beta)$ is said
to be
integrable if $\alpha$ and $\beta$ are closed.
\end{Def}

Note that a generalized Calabi-Yau structure is a Calabi-Yau structure
if and only if the metric is  definite positive and $K$ is a complex
structure. 

\begin{Rem}
The condition $d\alpha=d\beta=0$ implies the
integrability (in the Frobenius sense) of the distributions defined by
the almost complex structure or almost product structure
$K$. Therefore, according to the Newlander-Nirenberg theorem, it
implies its integrability. For instance, when $K$ is an almost complex
structure and $g$ is definite positive , the almost Calabi-Yau
structure $(g,\Omega,K,\alpha,\overline{\alpha})$ is integrable if and
only if $K$ is a complex structure and $\alpha$ is holomorphic.
\end{Rem}

\begin{Ex}
Each  nondegenerate Monge-Amp\`ere structure 
$(\Omega,\omega_0)$ defines the ge\-neralized almost Calabi-Yau structure
$(q_\omega,\Omega,K_\omega,\alpha,\beta)$ with 
$$
\omega=\frac{\omega_0}{\sqrt[4]{|\lambda(\omega_0)|}}.
$$
 
For instance, on $\mathbb{R}^6$, the generalized Calabi-Yau structure
associated with the equation
$$
\Delta(f)-\hess(f)=0
$$
is the canonical Calabi-Yau structure of $\mathbb{C}^3$
$$
\begin{cases}
g= -\underset{j=1}{\overset{3}{\sum}} dx_j.dx_j+ dy_j.dy_j&\\ 
K=\underset{j=1}{\overset{3}{\sum}} \frac{\partial}{\partial
  y_j}\otimes dx_j-\frac{\partial}{\partial x_j}\otimes dy_j&\\ 
\Omega=\underset{j=1}{\overset{3}{\sum}} dx_j\wedge dy_j&\\
\alpha=dz_1\wedge dz_2\wedge dz_3&\\
\beta=\overline{\alpha}
\end{cases}
$$
The generalized Calabi-Yau associated with the equation 
$$
\square(f)+\hess(f)=0
$$
is the pseudo Calabi-Yau structure  
$$
\begin{cases}
q=dx_1.dx_1-dx_2.dx_2+dx_3.dx_3 +dy_1.dy_1-dy_2.dy_2+dx_3.dx_3&\\
K=\frac{\partial}{\partial
  x_1}\otimes dy_1 -\frac{\partial}{\partial y_1}\otimes dx_1+
\frac{\partial}{\partial y_2}\otimes dx_2-\frac{\partial}{\partial x_2}\otimes
dy_2-\frac{\partial}{\partial y_3}\otimes
dx_3+\frac{\partial}{\partial x_3}\otimes dy_3&\\ 
\Omega=\underset{j=1}{\overset{3}{\sum}} dx_j\wedge dy_j&\\
\alpha=dz_1\wedge dz_2\wedge dz_3&\\
\beta=\overline{\alpha}
\end{cases}
$$  
The generalized Calabi-Yau structure associated with the equation
$$
\hess(f)=1 
$$
is the ``real'' Calabi-Yau structure
$$
\begin{cases}
g=\underset{j=1}{\overset{3}{\sum}} dx_j.dy_j&\\
K= \underset{j=1}{\overset{3}{\sum}} \frac{\partial}{\partial
  x_j}\otimes dx_j-\frac{\partial}{\partial y_j}\otimes dy_j&\\
\Omega=\underset{j=1}{\overset{3}{\sum}} dx_j\wedge dy_j&\\
\alpha=dx_1\wedge dx_2\wedge dx_3&\\
\beta= dy_1\wedge dy_2\wedge dy_3&\\
\end{cases}
$$
A manifold endowed with a ``real'' Calabi-Yau structure is the
analogue of a ``Monge-Amp\`ere manifold'' in the Kontsevich and
Soibelman sense (\cite{KS}). A Monge-Ampere manifold is an affine
riemannian manifold $(M,g)$ such that locally 
$$  
g=\sum_{i,j} \frac{\partial^2 F}{\partial x_i\partial x_j} dx_i.dx_j,
$$
$F$ being a smooth function satisfying 
$$
\dete\Big(\frac{\partial^2 F}{\partial x_i\partial x_j}\Big)=\text{
  constant}. 
$$
In the ``real'' Calabi-Yau case we have such a potential $F$:
$$
g=\sum_{i,j}\frac{\partial^2 F}{\partial x_i\partial y_j}dx_i.dy_j,
$$
and $\dete\Big(\frac{\partial^2 F}{\partial x_i\partial
  y_j}\Big)=f(x)g(y)$ (see \cite{Ba2} for more details).  
\end{Ex}

Let $(\Omega,\omega)$ be a Monge-Amp\`ere structure with
$\lambda(\omega)=\pm 1$. Since $d\omega=d\hat{\omega}=0$ if and only
if $d\alpha=d\beta=0$, we have the obvious proposition:

\begin{Prop}
A generalized almost Calabi-Yau structure associated with a
nondegenerate Monge-Ampere structure is integrable if and only this
Mon\-ge-Amp\`ere structure is closed.
\end{Prop}

\subsection{Nondegenerate Monge-Amp\`ere equations}

Let us come back now to the differential equation associated with a
nondegenerate Monge-Amp\`ere structure $(\Omega,\omega)$ on a $6$-dimensional
manifold $X$. It is natural to ask if this equation is locally
symplectically  equivalent to one of these:
$$
\begin{cases}     
\hess(f)=1&\\
\Delta(f)-\hess(f)=0&\\
\square(f)+\hess(f)=0&\\
\end{cases}
$$

According to table \ref{table1}, it will be the case if and only if
$(\Omega,\omega)$ is locally constant. The following theorem gives a
criterion using the generalized Calabi-Yau structure associated.

\begin{Theo}{\label{theo2}}
A Monge-Amp\`ere equation associated with a nondegenerate
Monge-Amp\`ere structure can be reduced by a symplectic change of
coordinates to one of the following equations
$$
\begin{cases}     
\hess(f)=1&\\
\Delta(f)-\hess(f)=0&\\
\square(f)+\hess(f)=0&\\
\end{cases}
$$
if and only if the generalized Calabi-Yau structure associated is
integrable and flat.
\end{Theo}

We refer to \cite{Ba2} for the proof. The idea is that the
integrability condition implies the existence of a ``generalized''
K\"ahler potential and the flat condition allows us to choose a
Darboux coordinates system in which this potential has a nice
expression.

Lychagin and Roubtsov have proved an equivalent theorem in
\cite{LR3} using technics of formal integrability. Theorem \ref{theo2} is more
restrictive since it only concerns nondegenerate Monge-Amp\`ere
equations but it is worth mentioning that its statement and its proof
are much more simple and that it has a nice geometric meaning.   

We sum up in  table \ref{table2} the correspondance between (pseudo)
Calabi-Yau structures and ellipic Monge-Amp\`ere structures.
\begin{table}[hbp!]
\begin{center}
\begin{tabular}{|c|c|}
\hline
almost (pseudo) CY& elliptic MA\\
\hline
(pseudo) CY& closed elliptic MA\\
\hline
flat (pseudo) CY& locally constant elliptic MA\\
\hline
\end{tabular}
\end{center}
\caption{(pseudo) Calabi-Yau and elliptic Monge-Amp\`ere
  structures}{\label{table2}} 
\end{table}

\begin{Ex}
There are very few explicit examples of Calabi-Yau metrics. One of
these is the Stenzel metric on $T^*S^n$ (see for instance
\cite{A},\cite{St}). This metric is not flat, therefore the special
lagrangian equation associated with is not the classical one. 

$T^*S^n=\Big\{(u,v)\in \mathbb{R}^{n+1}\times \mathbb{R}^{n+1}: \|u\|=1,
<u,v>=0\Big\}$ can be seen as the complex manifold  $Q^n=\Big\{z\in
C^{n+1}: z_1^2+\ldots + z_{n+1}^2=1\Big\}$ using the isomorphism
$$
\xi(x+iy)= (\frac{x}{\sqrt{1+\|y\|^2}},y).
$$
The holomorphic form is then 
$$
\alpha_z(Z_1,\ldots,Z_n)=det_{\mathbb{C}}(z,Z_1,\ldots,Z_n).
$$
and the K\"ahler form is $\Omega= i\partial \bar{\partial} \phi$ with
$\phi=f(\tau)$ where $\tau$ is the restriction to $Q^n$ of $|z_1|^2+
\ldots + |z_{n}^2$ and  $f$ is a solution of the ordinary differential
equation
$$
x(f')^n+ f''(f')^{n-1}(x^2-1)=c>0.
$$
To write the special lagrangian equation we have to find some Darboux
coordinates. Using the relations
$$
\begin{cases}
\underset{k=1}{\overset{4}{\sum}} u_kdu_k+v_kdv_k=0&\\
\underset{k=1}{\overset{4}{\sum}} u_kdv_k+v_kdu_k=0&\\
\end{cases}
$$
on $T^*S^3$, we see that on the chart $u_4\neq 0$,  
$$
\Omega=\sum_{k=1}^3 dw_k\wedge du_k
$$
with  
$$
w_k =2\frac{f'(2+2\|v\|^2)\sqrt{1+\|v\|^2}}{u_4}(u_kv_4-v_ku_4). 
$$
Denote by $\psi$ the map $(u,w)\mapsto (x+iy)$. The special lagrangian
equation on $T^*S^3$ is then
$$
(\psi\circ df)^*(Im(\alpha))=0.
$$ 
Note that it is difficult to  explicit this equation and it doesn't seem
possible to write it in a simple way.   
\end{Ex}

To resume, we have defined the generalized Calabi-Yau structures in order
to study an  equivalence problem  for Monge-Amp\`ere equations. The
author hopes that their construction is enough ``natural'' to have
a physical meaning. It is actually natural to ask if the generalized
Calabi-Yau manifolds have also mirror partners and if this mirror
conjecture can be formulated in terms of Monge-Amp\`ere equations.

\end{document}